\documentclass[12pt]{article}

\usepackage{amsmath,amsfonts,amssymb}

\newtheorem{theorem}{Theorem}[section]
\newtheorem{lemma}[theorem]{Lemma}

\newtheorem{proposition}[theorem]{Proposition}

\begin{document}

\title{\vspace*{-1.5cm}
Wong-Zakai Approximations of Backward Doubly Stochastic Doubly Backward Differential Equations}

\author{Ying Hu$^{1}$, Anis Matoussi$^{2}$ and
Tusheng Zhang$^{3}$}

\footnotetext[1]{\ Universit\'e{} de Rennes 1,
campus Beaulieu, 35042 Rennes Cedex, France}

\footnotetext[2]{\ Laboratoire Manceau de Math\'e{}matiques, Universit\'e{} du Maine,
Avenue Olivier Messiaen, 72 085 LE MANS Cedex, France}
\footnotetext[3]{\ School
of Mathematics, University of Manchester, Oxford Road, Manchester
M13 9PL, England, U.K. Email: tusheng.zhang@manchester.ac.uk}

\maketitle

\begin{abstract}
In this paper we obtain a Wong-Zakai approximation to solutions of backward doubly stochastic
differential equations.
\end{abstract}


\section{Framework and Introduction}
\setcounter{equation}{0}
Let $\{W_t,0\leq t\leq T\}$ and $\{B_t,0\leq t\leq T\}$ be two independent standard Brownian motions on a probability space $(\Omega, {\cal F}, P)$. Let ${\cal N}$ denote the class of $P$-null sets. For each $t\in [0,T]$, we define
$${\cal F}_t={\cal F}_t^W\vee {\cal F}_{t,T}^B,$$
where for any process $\{\eta_t\}$, ${\cal F}_{s, t}^{\eta}=\sigma\{\eta_r-\eta_s; s\leq r\leq t\}\vee {\cal N}$,${\cal F}_{t}^{\eta}={\cal F}_{0, t}^{\eta}$.
Let $f: R\times R\rightarrow R$ be a bounded measurable function satisfying
\vskip 0.3cm
{\bf (H.1)}
$$|f(y_1,z_1)-f(y_2,z_2)|\leq c (|y_1-y_2|+|z_1-z_2|)$$
Let $g\in C_b^1(R)$. For $n\geq 1$, define the linear interpolation $B^n$ of $B$ as
\begin{equation}\label{1.1}
B^n_t=B_{\frac{k+2}{2^n}}+2^n(t-\frac{k+1}{2^n})(B_{\frac{k+2}{2^n}}-B_{\frac{k+1}{2^n}}), \quad \mbox{for}\quad t\in [\frac{k}{2^n}, \frac{k+1}{2^n}].
\end{equation}
Let $\xi \in L^2(\Omega)$ be ${\cal F}_T$-measurable. Consider the following backward doubly stochastic differential equations(BDSDE):
\begin{eqnarray}\label{1.2}
Y_t&=& \xi+\int_t^Tf(Y_s,Z_s)ds+\int_t^Tg(Y_s)dB_s\nonumber\\
&&+\frac{1}{2}\int_t^Tgg^{\prime}(Y_s)ds-\int_t^TZ_sdW_s.
\end{eqnarray}
\begin{eqnarray}\label{1.3}
Y_t^n&=& \xi+\int_t^Tf(Y_s^n,Z_s^n)ds+\int_t^Tg(Y_s^n)dB_s^n\nonumber\\
&&-\int_t^TZ_s^ndW_s.
\end{eqnarray}
Here $dB_s$,$dB_s^n$ stand for the backward integrals. In the sequel, we will
use $\overrightarrow{ds}$ to indicate the backward integral against the Lebesgue measure.
\vskip 0.3cm
Backward doubly stochastic
differential equations was first studied by Pardoux and Peng in \cite{PP}. It is now a powerful tool to study stochastic partial differential equations with singular coefficients. Our purpose of this paper is to obtain the convergence of the Wong-Zakai approximation to the backward doubly stochastic differential equations, namely we will prove that $(Y^n,Z^n)$ converges to $(Y,Z)$ in $L^2$. The convergence of Wong-Zakai approximations to stochastic differential equations is now well known, see e.g. [IW]. Because of the nature of the BDSDEs, the integrand $Z^n, Z$ in the stochastic integral against Brownian motion are also part of the solutions. This makes the problem drastically different from the Wong-Zakai approximation for stochastic different equations. Another difficulty comes from the fact that the H\"o{}lder type estimate
$$ E[|Y_t^n-Y_s^n|^p]\leq C|t-s|^{\alpha}$$
is no longer available. We overcome this by carefully exploiting the independence of the two Brownian motions $B$ and $W$.
\vskip 0.4cm
The application of our results to stochastic partial differential equations will be discussed in a forthcoming paper. 
\section{Main results}
The following is an priori estimate for the family $\{(Y^n,Z^n), n\geq 1\}$.
\begin{proposition}
There exists a constant $C$ such that
\begin{equation}\label{2.1}
\sup_n\sup_{0\leq t\leq T}\{ E[(Y_t^n)^2]+E[\int_t^T(Z_s^n)^2ds]\}\leq C.
\end{equation}
\end{proposition}
\vskip 0.3cm
{\bf Proof}.
By Ito's formula, we have
\begin{eqnarray}\label{2.2}
&&(Y_t^n)^2+\int_t^T(Z_s^n)^2ds\nonumber\\
&=& (\xi)^2+2\int_t^TY_s^nf(Y_s^n,Z_s^n)ds+2\int_t^TY_s^ng(Y_s^n)dB_s^n\nonumber\\
&&-2\int_t^TY_s^nZ_s^ndW_s.
\end{eqnarray}
For $s\in [\frac{k}{2^n}, \frac{k+1}{2^n}]$, set $s^+=\frac{k+2}{2^n}$ and $s^-=\frac{k-1}{2^n}$.
In view of (H.1), it is easy to see that there exists $C_{1}>0$
such that
\begin{equation}\label{2.3}
2\int_t^TY_s^nf(Y_s^n,Z_s^n)ds\leq C_{1}\int_t^T(Y_s^n)^2ds+\frac{1}{4}\int_t^T(Z_s^n)^2ds
+C_1.
\end{equation}
Now, the third term on the right side of (\ref{2.2}) can be written as
\begin{eqnarray}\label{2.4}
2\int_t^TY_s^ng(Y_s^n)dB_s^n&=&2\int_t^TY_{s^+}^ng(Y_{s^+}^n)dB_s^n\nonumber\\
&+&2\int_t^T(Y_s^n-Y_{s^+}^n)g(Y_s^n)dB_s^n\nonumber\\
&+&2\int_t^TY_{s^+}^n(g(Y_s^n)-g(Y_{s^+}^n))dB_s^n\nonumber\\
&:=& I_1+I_2+I_3.
\end{eqnarray}
As a stochastic integral, we have $E[I_1]=0$. By the equation (\ref{1.3}) it follows that
\begin{eqnarray}\label{2.5}
I_2&=&2\int_t^T(\int_s^{s^+}f(Y^n_u,Z^n_u)du)g(Y_s^n)dB_s^n\nonumber\\
&+&2\int_t^T(\int_s^{s^+}g(Y^n_u)dB^n_u)g(Y_s^n)dB_s^n \nonumber\\
&-&2\int_t^T(\int_s^{s^+}Z_u^ndW_u)g(Y_s^n)dB_s^n\nonumber\\
&:=&I_{2.1}+I_{2.2}+I_{2.3}.
\end{eqnarray}
By the boundedness of $f$ and $g$, we have
\begin{eqnarray}\label{2.6}
E[I_{2.1}]&\leq &C\int_t^T(\int_s^{s^+}du)E[|\dot{B}_s^n|]ds\nonumber\\
&\leq & C (\frac{1}{2^n})^{\frac{1}{2}},
\end{eqnarray}
and
\begin{eqnarray}\label{2.7}
E[I_{2.2}]&\leq &CE[\int_t^T(\int_s^{s^+}|\dot{B}_u^n|du)|\dot{B}_s^n|ds]\nonumber\\
&\leq &C\int_t^Tds\int_s^{s^+}du E[|\dot{B}_u^n||\dot{B}_s^n|]\nonumber\\
&\leq &C\int_t^Tds\int_s^{s^+}du (E[|\dot{B}_u^n|^2])^{\frac{1}{2}}(E[|\dot{B}_s^n|^2])^{\frac{1}{2}}\nonumber\\
&\leq &C\int_t^Tds\int_s^{s^+}du (2^n)^{\frac{1}{2}}(2^n)^{\frac{1}{2}}\leq C,
\end{eqnarray}
where $C$ is a constant independent of $n$. For the term $I_{2.3}$, we have
\begin{eqnarray}\label{2.8}
E[I_{2.3}]&=&2\int_t^TE[(\int_s^{s^+}Z_u^ndW_u)g(Y_s^n)\dot{B}_s^n]ds\nonumber\\
&=&2\int_t^TE[g(Y_s^n)\dot{B}_s^nE[(\int_s^{s^+}Z_u^ndW_u)|{\cal F}_s]]ds  \nonumber\\
&= &0.
\end{eqnarray}
Putting together (\ref{2.6}),(\ref{2.7}),(\ref{2.8}) we get
\begin{equation}\label{2.9}
\sup_nE[I_2]\leq C,
\end{equation}
for some constant $C$.
To bound $I_3$ in (\ref{2.4}), we write
\begin{eqnarray}\label{2.10}
I_3&=&2\int_t^TY_{s^+}^n\int_0^1d\lambda g^{\prime}(Y_{s^+}^n+\lambda (Y_s^n-Y_{s^+}^n))(Y_s^n-Y_{s^+}^n)dB_s^n\nonumber\\
&=& 2\int_t^TY_{s^+}^n\int_0^1d\lambda g^{\prime}(Y_{s^+}^n+\lambda (Y_s^n-Y_{s^+}^n))(\int_s^{s^+}f(Y^n_u,Z^n_u)du)dB_s^n\nonumber\\
&+&2\int_t^TY_{s^+}^n\int_0^1d\lambda g^{\prime}(Y_{s^+}^n+\lambda (Y_s^n-Y_{s^+}^n))(\int_s^{s^+}g(Y^n_u)dB_u^n)dB_s^n\nonumber\\
&-&2\int_t^TY_{s^+}^n\int_0^1d\lambda g^{\prime}(Y_{s^+}^n+\lambda (Y_s^n-Y_{s^+}^n))(\int_s^{s^+}Z_u^ndW_u)dB_s^n\nonumber\\
&:=&I_{3.1}+I_{3.2}+I_{3.3}.
\end{eqnarray}
For $I_{3.1}$, we have
\begin{eqnarray}\label{2.11}
E[I_{3.1}]&\leq &C\int_t^TE[|Y_{s^+}^n|\frac{1}{2^n} |\dot{B}_s^n|]ds\nonumber\\
&\leq & C\int_t^TE[(Y_{s^+}^n)^2]ds+C\int_t^TE[(\frac{1}{2^n})^2 |\dot{B}_s^n|^2]ds\nonumber\\
&\leq & C\int_t^TE[(Y_{s^+}^n)^2]ds+C.
\end{eqnarray}
Similarly we have
\begin{eqnarray}\label{2.12}
E[I_{3.2}]&\leq &C\int_t^Tds\int_s^{s^+}du E[|Y_{s^+}^n||\dot{B}_u^n||\dot{B}_s^n|]ds\nonumber\\
&\leq & C\int_t^Tds\int_s^{s^+}du (E[(Y_{s^+}^n)^2])^{\frac{1}{2}}(E[|\dot{B}_u^n|^2 |\dot{B}_s^n|^2])^{\frac{1}{2}}ds\nonumber\\
&\leq & C\int_t^Tds\int_s^{s^+}du (E[(Y_{s^+}^n)^2])^{\frac{1}{2}}2^n\nonumber\\
&\leq & C\int_t^TE[(Y_{s^+}^n)^2]ds+C.
\end{eqnarray}
By virtue of the independence of  $Y_{s^+}^n$ and $\dot{B}_s^n$, we have
\begin{eqnarray}\label{2.13}
E[I_{3.3}]&\leq &CE[\int_t^T|Y_{s^+}^n||\dot{B}_s^n||\int_s^{s^+}Z_u^ndW_u|ds]
\nonumber\\
&\leq &C\int_t^T(E[(Y_{s^+}^n)^2|\dot{B}_s^n|^2])^{\frac{1}{2}}(E[|\int_s^{s^+}Z_u^ndW_u|^2])^{\frac{1}{2}}ds\nonumber\\
&= & C\int_t^T(E[(Y_{s^+}^n)^2])^{\frac{1}{2}}(E[|\dot{B}_s^n|^2])^{\frac{1}{2}}
(E[\int_s^{s^+}(Z_u^n)^2du])^{\frac{1}{2}}ds\nonumber\\
\nonumber\\
&= & C\int_t^T(E[(Y_{s^+}^n)^2])^{\frac{1}{2}}(2^n)^{\frac{1}{2}}(E[\int_s^{s^+}(Z_u^n)^2du])^{\frac{1}{2}}ds\nonumber\\
\nonumber\\
&\leq & \frac{1}{4}\int_t^T2^nE[\int_s^{s^+}(Z_u^n)^2du]ds+C_2\int_t^TE[(Y_{s^+}^n)^2]ds\nonumber\\
&=&\frac{1}{4}E[\int_t^T(Z_u^n)^2du 2^n(\int_{u^-}^{u}ds)]+C_2\int_t^TE[(Y_{s^+}^n)^2]ds \nonumber\\
&\leq &\frac{1}{4}E[\int_t^T(Z_u^n)^2du]+C_2\int_t^TE[(Y_{s^+}^n)^2]ds
\end{eqnarray}
(\ref{2.11})--(\ref{2.13}) imply that
\begin{equation}\label{2.14}
E[I_3]\leq \frac{1}{4}E[\int_t^T(Z_u^n)^2du]+C_2\int_t^TE[(Y_{s^+}^n)^2]ds+C
\end{equation}
It follows from (\ref{2.2}),(\ref{2.3}),(\ref{2.9}) and (\ref{2.14}) that
\begin{eqnarray}\label{2.15}
&&E[(Y_t^n)^2]+\frac{1}{2}E[\int_t^T(Z_s^n)^2ds]\nonumber\\
&\leq & E[(\xi)^2]+C\int_t^TE[(Y_{s^+}^n)^2]ds+C
\end{eqnarray}
Applying the Gronwall's inequality, we complete the proof of the Proposition.
\vskip 0.3cm
The above result can be strengthened as
\begin{proposition}
For any $p\geq 1$, there exists a constant $C$ such that
\begin{equation}\label{2.1-1}
\sup_n\{ E[\sup_{0\leq t\leq T}|Y_t^n|^p]+E[(\int_t^T(Z_s^n)^2ds)^p]\}\leq C.
\end{equation}
\end{proposition}

\begin{theorem}
\begin{equation}\label{2.16}
\lim_{n\rightarrow \infty}\sup_{0\leq t\leq T}\{E[(Y_t^n-Y_t)^2]+E[\int_t^T(Z_s^n-Z_s)^2ds]\}=0.
\end{equation}
\end{theorem}
{\bf Proof}.
By Ito's formula, we have
\begin{eqnarray}\label{2.17}
&&(Y_t^n-Y_t)^2+\int_t^T(Z_s^n-Z_s)^2ds\nonumber\\
&=& 2\int_t^T(Y_s^n-Y_s)(f(Y_s^n,Z_s^n)-f(Y_s,Z_s)) ds\nonumber\\
&+&2\int_t^T(Y_s^n-Y_s)g(Y_s^n)dB_s^n+\int_t^Tg^2(Y_s)ds\nonumber\\
&-&2\int_t^T(Y_s^n-Y_s)g(Y_s)dB_s-2\int_t^T(Y_s^n-Y_s)(Z_s^n-Z_s)dW_s\nonumber\\
&-&\int_t^T(Y_s^n-Y_s)gg^{\prime}(Y_s)ds\nonumber\\
&:=&I_1^n+I_2^n+I_3^n+I_4^n+I_5^n+I_6^n.
\end{eqnarray}
%Taking expectation in the above equation we get
%\begin{eqnarray}\label{2.18}
%&&E[(Y_t^n-Y_t)^2]+E\left [\int_t^T(Z_s^n-Z_s)^2ds\right ]\nonumber\\
%&=& 2E\left[\int_t^T(Y_s^n-Y_s)(f(Y_s^n,Z_s^n)-f(Y_s,Z_s)) ds\right]\nonumber\\
%&+&2E\left[\int_t^T(Y_s^n-Y_s)g(Y_s^n)dB_s^n\right]+E\left[\int_t^Tg^2(Y_s)ds\right]\nonumber\\
%&-&E\left [\int_t^T(Y_s^n-Y_s)gg^{\prime}(Y_s)ds\right ]\nonumber\\
%&:=& I_1^n+I_2^n+I_3^n+I_4^n.
%\end{eqnarray}

Crucially we need to bound the term $I_2^n$. We write
\begin{eqnarray}\label{2.18}
I_2^n&=&2\int_t^T(Y_s^n-Y_s)g(Y_s^n)dB_s^n\nonumber\\
&=&2\int_t^T[(Y_s^n-Y_s)-(Y_{s^+}^n-Y_{s^+})]g(Y_s^n)dB_s^n\nonumber\\
&+&2\int_t^T(Y_{s^+}^n-Y_{s^+})(g(Y_s^n)-g(Y_{s^+}^n))dB_s^n\nonumber\\
&+&2\int_t^T(Y_{s^+}^n-Y_{s^+})g(Y_{s^+}^n)dB_s^n\nonumber\\
&:=& A+B+C.
\end{eqnarray}
As a stochastic integral, we have $E[C]=0$. It is a long proof to establish the bounds for $E[A]$ and $E[B]$. We will split it into two lemmas for clarity.
 \begin{lemma}
 We have
 \begin{eqnarray}\label{2.19-1}
E[A]&\leq &C(\frac{1}{2^n})^{\frac{1}{2}-\delta}-2E[\int_t^Tg(Y_{s^+})g(Y^n_{s^+})ds]\nonumber\\
&&\quad\quad+E[\int_t^Tg^2(Y^n_{s^+})ds].
\end{eqnarray}
\end{lemma}
\vskip 0.3cm
{\bf Proof}. By the equations satisfied by $Y^n$ and $Y$ we have
\begin{eqnarray}\label{2.19}
A&=&2\int_t^T(\int_s^{s^+}g(Y^n_u)dB^n_u)g(Y_s^n)dB_s^n\nonumber\\
&-&2\int_t^T(\int_s^{s^+}g(Y_u)dB_u)g(Y_s^n)dB_s^n\nonumber\\
&+&2\int_t^T(\int_s^{s^+}[f(Y^n_u,Z_u^n)-f(Y_u,Z_u)]du)g(Y_s^n)dB_s^n\nonumber\\
&-&2\int_t^T(\int_s^{s^+}(Z_u^n-Z_u)dW_u)g(Y_s^n)dB_s^n\nonumber\\
&-&\int_t^T(\int_s^{s^+}gg^{\prime}(Y_u) du)g(Y_s^n)dB_s^n\nonumber\\
&:=&A_1+A_2+A_3+A_4+A_5.
\end{eqnarray}
Clearly,
\begin{equation}\label{2.20}
E[|A_5|]\leq C \frac{1}{2^n} \int_t^TE[|\dot{B}_s^n|]ds\leq C (\frac{1}{2^n})^{\frac{1}{2}},
\end{equation}
also
 \begin{equation}\label{2.20-1}
E[|A_3|]\leq C (\frac{1}{2^n})^{\frac{1}{2}}.
\end{equation}
By conditioning on ${\cal F}_s$, we find that
\begin{eqnarray}\label{2.21}
E[A_4]&=&-2\int_t^TE\left[(\int_s^{s^+}(Z_u^n-Z_u)dW_u)g(Y_s^n)\dot{B}_s^n\right]ds\nonumber\\
&=&-2\int_t^TE\left[E[(\int_s^{s^+}(Z_u^n-Z_u)dW_u)|{\cal F}_s]g(Y_s^n)\dot{B}_s^n\right]ds\nonumber\\
&=&0.
\end{eqnarray}
To bound $A_1$, we write it as
\begin{eqnarray}\label{2.22}
A_1&=&2\int_t^T(\int_s^{s^+}g(Y^n_u)dB^n_u)g(Y_s^n)dB_s^n\nonumber\\
&=&2\int_t^T[\int_s^{s^+}(g(Y^n_u)-g(Y^n_{s^+}))dB^n_u]g(Y_s^n)dB_s^n\nonumber\\
&+&2\int_t^Tg(Y^n_{s^+})(B_s^n-B_{s^+}^n)(g(Y^n_s)-g(Y^n_{s^+}))dB_s^n\nonumber\\
&+&2\int_t^Tg(Y^n_{s^+})(B_s^n-B_{s^+}^n)g(Y^n_{s^+})dB_s^n\nonumber\\
&:=&A_{11}+A_{12}+A_{13}.
\end{eqnarray}
Splitting the interval $[t,T]$ into subintervals $[\frac{k}{2^n}, \frac{k+1}{2^n}]$ we see that
\begin{eqnarray}\label{2.23}
A_{13}&=&2\sum_k\int_{\frac{k}{2^n}}^{\frac{k+1}{2^n}}g^2(Y^n_{\frac{k+2}{2^n}})2^n
(B_{\frac{k+2}{2^n}}-B_{\frac{k+3}{2^n}})(B_{\frac{k+2}{2^n}}-B_{\frac{k+1}{2^n}})
\overrightarrow{ds}\nonumber\\
&+&2\sum_k\int_{\frac{k}{2^n}}^{\frac{k+1}{2^n}}g^2(Y^n_{\frac{k+2}{2^n}})(2^n)^2
(s-\frac{k+1}{2^n})(B_{\frac{k+2}{2^n}}-B_{\frac{k+1}{2^n}})^2\overrightarrow{ds}\nonumber\\
&:=&A_{13,1}+A_{13,2}.
\end{eqnarray}
Conditioning on ${\cal F}_{\frac{k+2}{2^n}}$ we have
\begin{eqnarray}\label{2.24}
E[A_{13,1}]&=&-2\sum_kE[g^2(Y^n_{\frac{k+2}{2^n}})
(B_{\frac{k+2}{2^n}}-B_{\frac{k+3}{2^n}})E[(B_{\frac{k+2}{2^n}}-B_{\frac{k+1}{2^n}})|{\cal F}_{\frac{k+2}{2^n}}]] \nonumber\\
&=&0.
\end{eqnarray}
Integrating with respect to $s$, we get that
\begin{eqnarray}\label{2.25}
&&E[A_{13}]=E[A_{13,2}]\nonumber\\
&=&E[\sum_kg^2(Y^n_{\frac{k+2}{2^n}})
(B_{\frac{k+2}{2^n}}-B_{\frac{k+1}{2^n}})^2]\nonumber\\
&=&E[\sum_kg^2(Y^n_{\frac{k+2}{2^n}})\{
(B_{\frac{k+2}{2^n}}-B_{\frac{k+1}{2^n}})^2-\frac{1}{2^n}\} ]\nonumber\\
&&\quad\quad\quad +E[\int_t^Tg^2(Y^n_{s^+})ds]\nonumber\\
&=&E[\int_t^Tg^2(Y^n_{s^+})ds],
\end{eqnarray}
where the fact that the sequence $\{(B_{\frac{k+2}{2^n}}-B_{\frac{k+1}{2^n}})^2-\frac{1}{2^n}, k\geq 0\}$ is a martingale has been used.
For the term $A_{11}$ in (\ref{2.22}), we have
\begin{eqnarray}\label{2.26}
A_{11}&=&2\int_0^1d\lambda \int_t^T\int_s^{s^+}g^{\prime}(Y^n_{s^+}+\lambda (Y^n_u-Y^n_{s^+}))(Y^n_u-Y^n_{s^+})dB^n_ug(Y_s^n)dB_s^n\nonumber\\
&=&2\int_0^1d\lambda \int_t^T\int_s^{s^+}g^{\prime}(Y^n_{s^+}+\lambda (Y^n_u-Y^n_{s^+}))[\int_u^{s^+}f(Y^n_v,Z_v^n)dv]dB^n_ug(Y_s^n)dB_s^n\nonumber\\
&+&2\int_0^1d\lambda \int_t^T\int_s^{s^+}g^{\prime}(Y^n_{s^+}+\lambda (Y^n_u-Y^n_{s^+}))[\int_u^{s^+}g(Y^n_v)dB_v^n]dB^n_ug(Y_s^n)dB_s^n\nonumber\\
&-&2\int_0^1d\lambda \int_t^T\int_s^{s^+}g^{\prime}(Y^n_{s^+}+\lambda (Y^n_u-Y^n_{s^+}))[\int_u^{s^+}Z_v^ndW_v]dB^n_ug(Y_s^n)dB_s^n\nonumber\\
&:=&A_{11,1}+A_{11,2}+A_{11,3}.
\end{eqnarray}
The first two terms on the right can be bounded as follows.
\begin{eqnarray}\label{2.27}
E[A_{11,1}]
&\leq &C\frac{1}{2^n} \int_t^Tds\int_s^{s^+}E[|\dot{B}^n_u||\dot{B}^n_s|]du\nonumber\\
&\leq &C\frac{1}{2^n}(2^n)^{\frac{1}{2}}(2^n)^{\frac{1}{2}}\frac{1}{2^n}\leq C\frac{1}{2^n}.
\end{eqnarray}
\begin{eqnarray}\label{2.28}
E[A_{11,2}]
&\leq &C\int_t^Tds\int_s^{s^+}du\int_u^{s^+}dv E[|\dot{B}^n_u||\dot{B}^n_v||\dot{B}^n_s|]\nonumber\\
&\leq &C(2^n)^{\frac{3}{2}}(\frac{1}{2^n})^2\leq C(\frac{1}{2^n})^{\frac{1}{2}}.
\end{eqnarray}
The last term $A_{11,3}$ can be estimated as follows.
\begin{eqnarray}\label{2.29}
E[A_{11,3}]&\leq&CE[\int_t^Tds\int_s^{s^+}|\dot{B}^n_u||\dot{B}^n_s||\int_u^{s^+}Z_v^ndW_v|du]\nonumber\\
&\leq &C\int_t^Tds\int_s^{s^+}du (E[|\dot{B}^n_u|^2|\dot{B}^n_s|^2])^{\frac{1}{2}}(E[|\int_u^{s^+}Z_v^ndW_v|^2])^{\frac{1}{2}}
\nonumber\\
&\leq &C2^n\int_t^Tds\int_s^{s^+}du (E[\int_u^{s^+}(Z_v^n)^2dv])^{\frac{1}{2}}\nonumber\\
&\leq& C\int_t^Tds(E[\int_s^{s^+}(Z_v^n)^2dv])^{\frac{1}{2}}\nonumber\\
&\leq& C\left(\int_t^TdsE[\int_s^{s^+}(Z_v^n)^2dv]\right )^{\frac{1}{2}}\leq C\left(E[\int_t^T(Z_v^n)^2dv\int_{v^-}^{v}ds]\right )^{\frac{1}{2}}\nonumber\\
&\leq& C (\sup_n(E[\int_t^T(Z_v^n)^2dv)^{\frac{1}{2}})(\frac{1}{2^n})^{\frac{1}{2}}.
\end{eqnarray}
Putting together (\ref{2.27})--(\ref{2.29}) together we get
\begin{equation}\label{2.29-1}
E[A_{11}]\leq C(\frac{1}{2^n})^{\frac{1}{2}}.
\end{equation}
Similarly the term $A_{12}$ can be decomposed as
\begin{eqnarray}\label{2.30}
A_{12}&=&2\int_0^1d\lambda \int_t^Tg(Y_{s^+}^n)\int_s^{s^+}dB^n_ug^{\prime}(Y^n_{s^+}+\lambda (Y^n_s-Y^n_{s^+}))(Y^n_s-Y^n_{s^+})dB_s^n\nonumber\\
&=&2\int_0^1d\lambda \int_t^Tg(Y_{s^+}^n)\int_s^{s^+}dB^n_ug^{\prime}(Y^n_{s^+}+\lambda (Y^n_s-Y^n_{s^+}))[\int_s^{s^+}f(Y^n_v,Z_v^n)dv]dB_s^n\nonumber\\
&+&2\int_0^1d\lambda \int_t^Tg(Y_{s^+}^n)\int_s^{s^+}dB^n_ug^{\prime}(Y^n_{s^+}+\lambda (Y^n_s-Y^n_{s^+}))[\int_s^{s^+}g(Y^n_v)dB_v^n]dB_s^n \nonumber\\
&-&2\int_0^1d\lambda \int_t^Tg(Y_{s^+}^n)\int_s^{s^+}dB^n_ug^{\prime}(Y^n_{s^+}+\lambda (Y^n_s-Y^n_{s^+}))[\int_s^{s^+}Z_v^ndW_v]dB_s^n\nonumber\\
&:=&A_{12,1}+A_{12,2}+A_{12,3}.
\end{eqnarray}
Using the similar arguments as for (\ref{2.27}) and (\ref{2.28}) we can show that
\begin{equation}\label{2.31}
E[A_{12,j}]\leq C\frac{1}{2^n}, j=1,2,3.
\end{equation}
Hence,
\begin{equation}\label{2.33-1}
E[A_{12}]\leq C(\frac{1}{2^n})^{\frac{1}{2}}.
\end{equation}
Combining (\ref{2.25}),(\ref{2.29-1}) and (\ref{2.33-1}) we get
\begin{equation}\label{2.33-2}
E[A_{1}]\leq C(\frac{1}{2^n})^{\frac{1}{2}}+E[\int_t^Tg^2(Y^n_{s^+})ds].
\end{equation}
Now we turn to $A_2$ which can be written as
\begin{eqnarray}\label{2.34}
A_2&=&-2\int_t^T(\int_s^{s^+}g(Y_u)dB_u)g(Y_s^n)dB_s^n\nonumber\\
&=&-2\int_t^T[\int_s^{s^+}(g(Y_u)-g(Y_{s^+}))dB_u]g(Y_s^n)dB_s^n\nonumber\\
&-&2\int_t^Tg(Y_{s^+})(B_s-B_{s^+})(g(Y^n_s)-g(Y^n_{s^+}))dB_s^n\nonumber\\
&-&2\int_t^Tg(Y_{s^+})(B_s-B_{s^+})g(Y^n_{s^+})dB_s^n\nonumber\\
&:=&A_{21}+A_{22}+A_{23}.
\end{eqnarray}
Splitting the interval $[t,T]$ into subintervals $[\frac{k}{2^n}, \frac{k+1}{2^n}]$,
\begin{eqnarray}\label{2.35}
A_{23}&=&-2\sum_k\int_{\frac{k}{2^n}}^{\frac{k+1}{2^n}}g(Y_{\frac{k+2}{2^n}})
g(Y^n_{\frac{k+2}{2^n}})
(B_{s}-B_{\frac{k+2}{2^n}})2^n(B_{\frac{k+2}{2^n}}-B_{\frac{k+1}{2^n}})\overrightarrow{ds}\nonumber\\
&=&-2\sum_k\int_{\frac{k}{2^n}}^{\frac{k+1}{2^n}}g(Y_{\frac{k+2}{2^n}})
g(Y^n_{\frac{k+2}{2^n}})
(B_{s}-B_{\frac{k+1}{2^n}})2^n(B_{\frac{k+2}{2^n}}-B_{\frac{k+1}{2^n}})\overrightarrow{ds}\nonumber\\
&-&2\sum_kg(Y_{\frac{k+2}{2^n}})
g(Y^n_{\frac{k+2}{2^n}})(B_{\frac{k+2}{2^n}}-B_{\frac{k+1}{2^n}})^2\nonumber\\
&:=&A_{23,1}+A_{23,2}.
\end{eqnarray}
Using the independence of the increments of $B$ and conditioning on ${\cal F}_{\frac{k+2}{2^n}}$ it follows that
\begin{eqnarray}\label{2.36}
&&E[A_{23,1}]\nonumber\\
&=&-2\sum_k\int_{\frac{k}{2^n}}^{\frac{k+1}{2^n}}2^n \overrightarrow{ds} E[g(Y_{\frac{k+2}{2^n}})g(Y^n_{\frac{k+2}{2^n}})E[(B_{s}-B_{\frac{k+1}{2^n}})(B_{\frac{k+2}{2^n}}-B_{\frac{k+1}{2^n}})|{\cal F}_{\frac{k+2}{2^n}}]] \nonumber\\
&=&-2\sum_k\int_{\frac{k}{2^n}}^{\frac{k+1}{2^n}}2^n \overrightarrow{ds} E[g(Y_{\frac{k+2}{2^n}})g(Y^n_{\frac{k+2}{2^n}})E[(B_{s}-B_{\frac{k+1}{2^n}})
(B_{\frac{k+2}{2^n}}-B_{\frac{k+1}{2^n}})]] \nonumber\\
&=&0.
\end{eqnarray}
Now,
\begin{eqnarray}\label{2.37}
A_{23,2}&=&-2\sum_kg(Y_{\frac{k+2}{2^n}})
g(Y^n_{\frac{k+2}{2^n}})\{(B_{\frac{k+2}{2^n}}-B_{\frac{k+1}{2^n}})^2-\frac{1}{2^n}\}\nonumber\\
&-&2\int_t^Tg(Y_{s^+})g(Y^n_{s^+})ds
\end{eqnarray}
By conditioning on ${\cal F}_{\frac{k+2}{2^n}}$ we see that the expectation of the first term of the above equation vanishes. Hence,
\begin{eqnarray}\label{2.38}
&&E[A_{23}]=E[A_{23,2}]=-2E[\int_t^Tg(Y_{s^+})g(Y^n_{s^+})ds].
\end{eqnarray}
For the term $A_{21}$ we have
\begin{eqnarray}\label{2.39}
A_{21}&=&-2\int_0^1d\lambda \int_t^T\int_s^{s^+}g^{\prime}(Y_{s^+}+\lambda (Y_u-Y_{s^+}))(Y_u-Y_{s^+})dB_ug(Y_s^n)dB_s^n\nonumber\\
&=&-2\int_0^1d\lambda \int_t^T\int_s^{s^+}g^{\prime}(Y_{s^+}+\lambda (Y_u-Y_{s^+}))[\int_u^{s^+}f(Y_v,Z_v)dv]dB_ug(Y_s^n)dB_s^n\nonumber\\
&-&2\int_0^1d\lambda \int_t^T\int_s^{s^+}g^{\prime}(Y_{s^+}+\lambda (Y_u-Y_{s^+}))[\int_u^{s^+}g(Y_v)dB_v]dB_ug(Y_s^n)dB_s^n\nonumber\\
&+&2\int_0^1d\lambda \int_t^T\int_s^{s^+}g^{\prime}(Y_{s^+}+\lambda (Y_u-Y_{s^+}))[\int_u^{s^+}Z_vdW_v]dB_ug(Y_s^n)dB_s^n\nonumber\\
&:=&A_{21,1}+A_{21,2}+A_{21,3}.
\end{eqnarray}
We will estimate each of the terms on the right. First we have
\begin{eqnarray}\label{2.40}
E[A_{21,1}]&\leq &C\int_0^1d\lambda \int_t^TE\left [ |\int_s^{s^+}g^{\prime}(Y_{s^+}+\lambda (Y_u-Y_{s^+}))[\int_u^{s^+}f(Y_v,Z_v)dv]dB_u| |\dot{B}_s^n|\right]ds\nonumber\\
&\leq &C\int_0^1d\lambda \int_t^T(E[ |\int_s^{s^+}g^{\prime}(Y_{s^+}+\lambda (Y_u-Y_{s^+}))\nonumber\\
&&\quad\quad\times [\int_u^{s^+}f(Y_v,Z_v)dv]dB_u|^2])^{\frac{1}{2}}(E[|\dot{B}_s^n|^2])^{\frac{1}{2}}ds\nonumber\\
&\leq &C(2^n)^{\frac{1}{2}}\int_0^1d\lambda \int_t^T(E[ \int_s^{s^+}g^{\prime}(Y_{s^+}+\lambda (Y_u-Y_{s^+}))^2(\int_u^{s^+}f(Y_v,Z_v)dv)^2du])^{\frac{1}{2}}ds\nonumber\\
&\leq& C(2^n)^{\frac{1}{2}} \int_t^T[(s^+-s)(s^+-s)^2 ]^{\frac{1}{2}}ds\nonumber\\
&\leq& C(\frac{1}{2^n}).
\end{eqnarray}
Similarly,
\begin{eqnarray}\label{2.41}
&&E[A_{21,2}]\nonumber\\
&\leq &C\int_0^1d\lambda \int_t^TE\left [ |\int_s^{s^+}g^{\prime}(Y_{s^+}+\lambda (Y_u-Y_{s^+}))[\int_u^{s^+}g(Y_v)dB_v]dB_u| |\dot{B}_s^n|\right]ds\nonumber\\
&\leq &C\int_0^1d\lambda \int_t^T(E[ |\int_s^{s^+}g^{\prime}(Y_{s^+}+\lambda (Y_u-Y_{s^+}))[\int_u^{s^+}g(Y_v)dB_v]dB_u|^2])^{\frac{1}{2}}(E[|\dot{B}_s^n|^2])^{\frac{1}{2}}ds\nonumber\\
&\leq &C(2^n)^{\frac{1}{2}}\int_0^1d\lambda \int_t^T(E[ \int_s^{s^+}g^{\prime}(Y_{s^+}+\lambda (Y_u-Y_{s^+}))^2(\int_u^{s^+}g(Y_v)dB_v)^2du])^{\frac{1}{2}}ds\nonumber\\
&\leq&C(2^n)^{\frac{1}{2}}\int_t^T(\int_s^{s^+}E[(\int_u^{s^+}g(Y_v)dB_v)^2]du
)^{\frac{1}{2}}ds\nonumber\\
&\leq&C(2^n)^{\frac{1}{2}}\int_t^T(\int_s^{s^+}du)^{\frac{1}{2}}
(E[\int_s^{s^+}g^2(Y_v)dv)])^{\frac{1}{2}}ds\nonumber\\
&\leq& C(\frac{1}{2^n})^{\frac{1}{2}}.
\end{eqnarray}
By H\"o{}lder inequality and Ito isometry, we have
\begin{eqnarray}\label{2.42}
&&E[A_{21,3}]\nonumber\\
&\leq &C\int_0^1d\lambda \int_t^TE\left [ |\int_s^{s^+}g^{\prime}(Y_{s^+}+\lambda (Y_u-Y_{s^+}))[\int_u^{s^+}Z_vdW_v]dB_u| |\dot{B}_s^n|\right]ds\nonumber\\
&\leq &C\int_0^1d\lambda \int_t^T(E[ |\int_s^{s^+}g^{\prime}(Y_{s^+}+\lambda (Y_u-Y_{s^+}))[\int_u^{s^+}Z_vdW_v]dB_u|^2])^{\frac{1}{2}}(E[|\dot{B}_s^n|^2])^{\frac{1}{2}}ds\nonumber\\
&\leq &C(2^n)^{\frac{1}{2}}\int_0^1d\lambda \int_t^T(E[ \int_s^{s^+}g^{\prime}(Y_{s^+}+\lambda (Y_u-Y_{s^+}))^2(\int_u^{s^+}Z_vdW_v)^2du])^{\frac{1}{2}}ds\nonumber\\
&\leq&C(2^n)^{\frac{1}{2}}\int_t^T(\int_s^{s^+}E[(\int_u^{s^+}Z_vdW_v)^2]du
)^{\frac{1}{2}}ds\nonumber\\
&\leq&C(2^n)^{\frac{1}{2}}\int_t^T(\int_s^{s^+}du)^{\frac{1}{2}}
(E[\int_s^{s^+}Z_v^2dv])^{\frac{1}{2}}ds\nonumber\\
&\leq&C(\int_t^T(E[\int_s^{s^+}Z_v^2dv]ds )^{\frac{1}{2}}\nonumber\\
&\leq&C(E[\int_t^TZ_v^2dv\int_{v^-}^{v}ds])^{\frac{1}{2}}\nonumber\\
&\leq& C(\frac{1}{2^n})^{\frac{1}{2}}.
\end{eqnarray}
It follows from (\ref{2.40}), (\ref{2.41}) and (\ref{2.42})that
\begin{equation}\label{2.42-1}
E[A_{21}]\leq C(\frac{1}{2^n})^{\frac{1}{2}}.
\end{equation}
Let us turn to the term $A_{22}$. We have
\begin{eqnarray}\label{2.42-2}
A_{22}&\leq &C\int_t^T|B_s-B_{s^+}||Y_s^n-Y_{s^+}^n| |\dot{B}_s^n|ds\nonumber\\
&\leq &C\int_t^T|B_s-B_{s^+}||\int_s^{s^+}f(Y_u^n,Z_u^n)du| |\dot{B}_s^n|ds\nonumber\\
&+ &C\int_t^T|B_s-B_{s^+}||\int_s^{s^+}g(Y_u^n)dB_u^n| |\dot{B}_s^n|ds\nonumber\\
&+ &C\int_t^T|B_s-B_{s^+}||\int_s^{s^+}Z_u^ndW_u| |\dot{B}_s^n|ds\nonumber\\
&:=& A_{22,1}+A_{22,2}+A_{22,3}.
\end{eqnarray}
Now,
\begin{eqnarray}\label{2.42-3}
E[A_{22,1}]
&\leq &C\frac{1}{2^n}\int_t^T(E[|B_s-B_{s^+}|^2])^{\frac{1}{2}} (E[|\dot{B}_s^n|^2])^{\frac{1}{2}}ds\nonumber\\
&\leq&C\frac{1}{2^n},
\end{eqnarray}
and
\begin{eqnarray}\label{2.42-4}
E[A_{22,2}]
&\leq &C\frac{1}{2^n}\int_t^T(E[|B_s-B_{s^+}|^2])^{\frac{1}{2}} (E[\sup_u|\dot{B}_u^n|^4])^{\frac{1}{2}}ds\nonumber\\
&\leq &C(\frac{1}{2^n})^{\frac{1}{2}}(2^n)^2\int_t^T (E[\sup_{|r-v|\leq \frac{1}{2^n}}|B_r-B_v|^4])^{\frac{1}{2}}ds\nonumber\\
&\leq&C(\frac{1}{2^n})^{\frac{1}{2}-\delta}.
\end{eqnarray}
By H\"o{}lder inequality,
\begin{eqnarray}\label{2.42-5}
E[A_{22,3}]
&\leq &C\int_t^T(E[|B_s-B_{s^+}|^4])^{\frac{1}{4}} (E[|\dot{B}_s^n|^4])^{\frac{1}{4}} (E[|\int_s^{s^+}Z_u^ndW_u|^2])^{\frac{1}{2}}ds\nonumber\\
&\leq&C(\int_t^T E[\int_s^{s^+}(Z_u^n)^2du]ds)^{\frac{1}{2}}\nonumber\\
&\leq&C\frac{1}{2^n}.
\end{eqnarray}
It follows from (\ref{2.42-2})--(\ref{2.42-5}) that
\begin{equation}\label{2.42-6}
E[A_{22}]\leq C(\frac{1}{2^n})^{\frac{1}{2}-\delta}.
\end{equation}
Collecting (\ref{2.38}),(\ref{2.42-1}) and (\ref{2.42-2}) we obtain
\begin{equation}\label{2.42-7}
E[A_{2}]\leq C(\frac{1}{2^n})^{\frac{1}{2}-\delta}-2E[\int_t^Tg(Y_{s^+})g(Y^n_{s^+})ds].
\end{equation}
The proof  is now completed by putting (\ref{2.20}),(\ref{2.20-1}),(\ref{2.21}),(\ref{2.33-2})
and (\ref{2.42-7}) together.
\begin{lemma}
 We have
 \begin{eqnarray}\label{2.42-8}
E[B]&\leq &C(\frac{1}{2^n})^{\frac{1}{2}-\delta}+E[\int_t^T(Y^n_{s^+}-Y_{s^+})gg^{\prime}(Y^n_{s^+})ds].
\end{eqnarray}
\end{lemma}
{\bf Proof}. Write
\begin{eqnarray}\label{2.43}
B&=&2\int_t^T(Y_{s^+}^n-Y_{s^+})(g(Y_s^n)-g(Y_{s^+}^n))dB_s^n\nonumber\\
&=&2\int_0^1d\lambda \int_t^T(Y_{s^+}^n-Y_{s^+})g^{\prime}(Y^n_{s^+}+\lambda (Y^n_s-Y^n_{s^+}))(Y^n_s-Y^n_{s^+})dB_s^n\nonumber\\
&=&2\int_0^1d\lambda \int_t^T(Y_{s^+}^n-Y_{s^+})g^{\prime}(Y^n_{s^+}+\lambda (Y^n_s-Y^n_{s^+}))[\int_s^{s^+}f(Y^n_v,Z_v^n)dv]dB_s^n\nonumber\\
&+&2\int_0^1d\lambda \int_t^T(Y_{s^+}^n-Y_{s^+})g^{\prime}(Y^n_{s^+}+\lambda (Y^n_s-Y^n_{s^+}))[\int_s^{s^+}g(Y^n_v)dB_v^n]dB_s^n\nonumber\\
&-&2\int_0^1d\lambda \int_t^T(Y_{s^+}^n-Y_{s^+})g^{\prime}(Y^n_{s^+}+\lambda (Y^n_s-Y^n_{s^+}))[\int_s^{s^+}Z_v^ndW_v]dB_s^n\nonumber\\
&:=&B_1+B_2+B_3.
\end{eqnarray}
By Proposition 2.2, we have
\begin{eqnarray}\label{2.44}
E[B_1]
&\leq &C\frac{1}{2^n} E[ \int_t^T|Y_{s^+}^n-Y_{s^+}||\dot{B}_s^n|ds]\nonumber\\
&\leq &C\frac{1}{2^n} \int_t^T(E[|Y_{s^+}^n-Y_{s^+}|^2])^{\frac{1}{2}}(E[|\dot{B}_s^n|^2])^{\frac{1}{2}}ds\nonumber\\
&\leq& C(\frac{1}{2^n})^{\frac{1}{2}}.
\end{eqnarray}
$B_2$ is further written as follows.
\begin{eqnarray}\label{2.45}
B_2&=&2\int_0^1d\lambda \int_t^T(Y_{s^+}^n-Y_{s^+})[g^{\prime}(Y^n_{s^+}+\lambda (Y^n_s-Y^n_{s^+}))-g^{\prime}(Y^n_{s^+})][\int_s^{s^+}g(Y^n_u)dB_u^n]dB_s^n\nonumber\\
&+&2\int_0^1d\lambda \int_t^T(Y_{s^+}^n-Y_{s^+})g^{\prime}(Y^n_{s^+})
[\int_s^{s^+}(g(Y^n_u)-g(Y^n_{s^+}))dB_u^n]dB_s^n\nonumber\\
&+&2\int_0^1d\lambda \int_t^T(Y_{s^+}^n-Y_{s^+})g^{\prime}(Y^n_{s^+})g(Y^n_{s^+})
(B_s^n-B_{s^+}^n)dB_s^n\nonumber\\
&:=&B_{21}+B_{22}+B_{23}.
\end{eqnarray}
By the Lipschitz continuity of $g^{\prime}$, it follows that
\begin{eqnarray}\label{2.46}
B_{21}&\leq &C \int_t^T|Y_{s^+}^n-Y_{s^+}|| Y^n_s-Y^n_{s^+}||\int_s^{s^+}g(Y^n_u)dB_u^n||\dot{B}_s^n|ds\nonumber\\
&\leq &C \int_t^T|Y_{s^+}^n-Y_{s^+}||\int_s^{s^+}f(Y^n_v,Z_v^n)dv|
|\int_s^{s^+}g(Y^n_u)dB_u^n||\dot{B}_s^n|ds\nonumber\\
&+ &C \int_t^T|Y_{s^+}^n-Y_{s^+}||\int_s^{s^+}g(Y^n_u)dB_u^n|^2|\dot{B}_s^n|ds\nonumber\\
&+ &C \int_t^T|Y_{s^+}^n-Y_{s^+}||\int_s^{s^+}Z_v^ndW_v|
|\int_s^{s^+}g(Y^n_u)dB_u^n||\dot{B}_s^n|ds\nonumber\\
&:=&B_{21,1}+B_{21,2}+B_{21,3}.
\end{eqnarray}
The following two inequalities will be used frequently in sequel.
\begin{equation}\label{2.47}
\sup_u|\dot{B}_u^n|\leq 2^n \sup_{|r-s|\leq \frac{1}{2^n}}|B_r-B_s|
\end{equation}
For any $\delta>0$ and $p\geq 1$, there exists a constant $C_{p,\delta}$ such that
\begin{equation}\label{2.48}
E[\sup_{|r-s|\leq \frac{1}{2^n}}|B_r-B_s|^p]\leq C_{p,\delta} (\frac{1}{2^n})^{\frac{p}{2}-\delta}.
\end{equation}
By H\"o{}lder's inequality and (\ref{2.47}), (\ref{2.48}), we have
\begin{eqnarray}\label{2.49}
E[B_{21,1}]&\leq &C \frac{1}{2^n}E[\int_t^T|Y_{s^+}^n-Y_{s^+}||\int_s^{s^+}|\dot{B}_u^n|du||\dot{B}_s^n|ds]\nonumber\\
&\leq &C \frac{1}{2^n}(2^n)^2 E[\int_t^T|Y_{s^+}^n-Y_{s^+}||\int_s^{s^+}du|\sup_{|r-s|\leq \frac{1}{2^n}}|B_r-B_s|^2 ds]\nonumber\\
&\leq &C E[\int_t^T|Y_{s^+}^n-Y_{s^+}|\sup_{|r-s|\leq \frac{1}{2^n}}|B_r-B_s|^2 ds]\nonumber\\
&\leq &C \int_t^T(E[|Y_{s^+}^n-Y_{s^+}|^2])^{\frac{1}{2}}  (E[\sup_{|r-s|\leq \frac{1}{2^n}}|B_r-B_s|^4])^{\frac{1}{2}}  ds\nonumber\\
&\leq& C(\frac{1}{2^n})^{1-\delta}.
\end{eqnarray}
Similarly, in view of (\ref{2.48}), we have
\begin{eqnarray}\label{2.50}
E[B_{21,2}]&\leq &CE[\int_t^T|Y_{s^+}^n-Y_{s^+}||\int_s^{s^+}du|^2\sup_s|\dot{B}_s^n|^3 ds]\nonumber\\
&\leq &C (\frac{1}{2^n})^2(2^n)^3 E[\int_t^T|Y_{s^+}^n-Y_{s^+}|\sup_{|r-s|\leq \frac{1}{2^n}}|B_r-B_s|^3 ds]\nonumber\\
&\leq &C 2^n\int_t^T(E[|Y_{s^+}^n-Y_{s^+}|^2])^{\frac{1}{2}}  (E[\sup_{|r-s|\leq \frac{1}{2^n}}|B_r-B_s|^6])^{\frac{1}{2}}  ds\nonumber\\
&\leq& C(\frac{1}{2^n})^{\frac{1}{2}-\delta},
\end{eqnarray}
and
\begin{eqnarray}\label{2.51}
E[B_{21,3}]&\leq &CE[\int_t^T|Y_{s^+}^n-Y_{s^+}||\int_s^{s^+}Z_v^ndW_v|(\int_s^{s^+} |\dot{B}_u^n|du)|\dot{B}_s^n|ds]\nonumber\\
&\leq &C \frac{1}{2^n}(2^n)^2 E[\int_t^T|Y_{s^+}^n-Y_{s^+}||\int_s^{s^+}Z_v^ndW_v|\sup_{|r-s|\leq \frac{1}{2^n}}|B_r-B_s|^2 ds]\nonumber\\
&\leq &C 2^n\int_t^T(E[|Y_{s^+}^n-Y_{s^+}|^4])^{\frac{1}{4}}(E[|\int_s^{s^+}Z_v^ndW_v|^2])^{\frac{1}{2}}  (E[\sup_{|r-s|\leq \frac{1}{2^n}}|B_r-B_s|^8])^{\frac{1}{4}}  ds\nonumber\\
&\leq &C (\frac{1}{2^n})^{1-\delta} 2^n\int_t^T(E[|Y_{s^+}^n-Y_{s^+}|^4])^{\frac{1}{4}}(E[\int_s^{s^+}(Z_v^n)^2dv])^{\frac{1}{2}} ds\nonumber\\
&\leq &C (\frac{1}{2^n})^{1-\delta} 2^n(\int_t^T(E[|Y_{s^+}^n-Y_{s^+}|^4])^{\frac{1}{2}}ds)^{\frac{1}{2}}(\int_t^T E[\int_s^{s^+}(Z_v^n)^2dv]ds )^{\frac{1}{2}} \nonumber\\
&\leq &C (\frac{1}{2^n})^{1-\delta} 2^n(E[\int_t^T(Z_v^n)^2dv \int_{v^-}^{v}ds])^{\frac{1}{2}} \nonumber\\
&\leq& C(\frac{1}{2^n})^{\frac{1}{2}-\delta},
\end{eqnarray}
where the a priori estimate (\ref{2.1-1}) has been used. (\ref{2.49})--(\ref{2.49})yields
\begin{equation}\label{2.51-1}
E[B_{21}]\leq C(\frac{1}{2^n})^{\frac{1}{2}-\delta}.
\end{equation}

By the Lipschitz continuity of $g$, we have
\begin{eqnarray}\label{2.52}
B_{22}&\leq &C \int_t^T|Y_{s^+}^n-Y_{s^+}|(\int_s^{s^+}| Y^n_u-Y^n_{s^+}|\dot{B}_u^n|du)|\dot{B}_s^n|ds\nonumber\\
&\leq &C \int_t^T|Y_{s^+}^n-Y_{s^+}|(\int_s^{s^+}|\int_u^{s^+}f(Y^n_v,Z_v^n)dv|
|\dot{B}_u^n|du)|\dot{B}_s^n|ds\nonumber\\
&\leq &C \int_t^T|Y_{s^+}^n-Y_{s^+}|(\int_s^{s^+}|\int_u^{s^+}g(Y^n_v)dB_v^n|
|\dot{B}_u^n|du)|\dot{B}_s^n|ds\nonumber\\
&\leq &C \int_t^T|Y_{s^+}^n-Y_{s^+}|(\int_s^{s^+}|\int_u^{s^+}Z_v^ndW_v|
|\dot{B}_u^n|du)|\dot{B}_s^n|ds\nonumber\\
&:=&B_{22,1}+B_{22,2}+B_{22,3}.
\end{eqnarray}
By the similar arguments as above, we have
\begin{equation}\label{2.53}
E[B_{22,j}]\leq C(\frac{1}{2^n})^{1-\delta}, j=1,2,3.
\end{equation}
Hence,
\begin{equation}\label{2.56}
E[B_{22}]\leq C(\frac{1}{2^n})^{\frac{1}{2}-\delta}.
\end{equation}
For the term $B_{23}$, we split  the interval $[t,T]$ into subintervals $[\frac{k}{2^n}, \frac{k+1}{2^n}]$ to get
\begin{eqnarray}\label{2.57}
&&B_{23}\nonumber\\
&=&2\sum_k\int_{\frac{k}{2^n}}^{\frac{k+1}{2^n}}(Y^n_{\frac{k+2}{2^n}}-Y_{\frac{k+2}{2^n}})
gg^{\prime}(Y^n_{\frac{k+2}{2^n}})2^n
(B_{\frac{k+2}{2^n}}-B_{\frac{k+3}{2^n}})(B_{\frac{k+2}{2^n}}-B_{\frac{k+1}{2^n}})\overrightarrow{ds}\nonumber\\
&+&2\sum_k\int_{\frac{k}{2^n}}^{\frac{k+1}{2^n}}(Y^n_{\frac{k+2}{2^n}}-Y_{\frac{k+2}{2^n}})
gg^{\prime}(Y^n_{\frac{k+2}{2^n}})(2^n)^2
(s-\frac{k+1}{2^n})(B_{\frac{k+2}{2^n}}-B_{\frac{k+1}{2^n}})^2\overrightarrow{ds}\nonumber\\
&=&-2\sum_k(Y^n_{\frac{k+2}{2^n}}-Y_{\frac{k+2}{2^n}})
gg^{\prime}(Y^n_{\frac{k+2}{2^n}})
(B_{\frac{k+2}{2^n}}-B_{\frac{k+3}{2^n}})(B_{\frac{k+2}{2^n}}-B_{\frac{k+1}{2^n}})\nonumber\\
&+&\sum_k(Y^n_{\frac{k+2}{2^n}}-Y_{\frac{k+2}{2^n}})
gg^{\prime}(Y^n_{\frac{k+2}{2^n}})\{(B_{\frac{k+2}{2^n}}-B_{\frac{k+1}{2^n}})^2-\frac{1}{2^n}\}
\nonumber\\
&+&\int_t^T(Y^n_{s^+}-Y_{s^+})gg^{\prime}(Y^n_{s^+})ds
\end{eqnarray}
Conditioning on ${\cal F}_{\frac{k+2}{2^n}}$ it is easy to see that the expectation of the first two terms on the right vanishes. Hence,
\begin{equation}\label{2.58}
E[B_{23}]=E[\int_t^T(Y^n_{s^+}-Y_{s^+})gg^{\prime}(Y^n_{s^+})ds].
\end{equation}
Collect the terms in (\ref{2.51-1}),(\ref{2.56}) and (\ref{2.58}) to obtain
\begin{equation}\label{2.58-1}
E[B_{2}]\leq C(\frac{1}{2^n})^{\frac{1}{2}-\delta}+E[\int_t^T(Y^n_{s^+}-Y_{s^+})gg^{\prime}(Y^n_{s^+})ds].
\end{equation}
Now we turn to the term $B_3$ in (\ref{2.43}). We further split it as
\begin{eqnarray}\label{2.59}
B_3&=&-2\{\int_0^1d\lambda \int_t^T[(Y_{s^+}^n-Y_{s^+})g^{\prime}(Y^n_{s^+}+\lambda (Y^n_s-Y^n_{s^+}))\nonumber\\
&&\quad -(Y_{s}^n-Y_{s})g^{\prime}(Y^n_{s})](\int_s^{s^+}Z_v^ndW_v)dB_s^n\}\nonumber\\
&&-2\int_t^T(Y_{s}^n-Y_{s})g^{\prime}(Y^n_{s})(\int_s^{s^+}Z_v^ndW_v)dB_s^n\nonumber\\
&:=&B_{31}+B_{32}.
\end{eqnarray}
First we notice that
\begin{eqnarray}\label{2.60}
E[B_{32}]
&=&-2\int_t^TE[(Y_{s}^n-Y_{s})g^{\prime}(Y^n_{s})\dot{B}_s^nE[(\int_s^{s^+}Z_v^ndW_v)|{\cal F}_s]]ds\nonumber\\
&=&0.
\end{eqnarray}
By the Lipschitz continuity of $g^{\prime}$, we have
\begin{eqnarray}\label{2.61}
B_{31}&\leq &C \int_t^T|Y_{s^+}^n-Y_s^n||\int_s^{s^+}Z_v^ndW_v||\dot{B}_s^n|ds\nonumber\\
&+&C \int_t^T|Y_{s^+}-Y_s||\int_s^{s^+}Z_v^ndW_v||\dot{B}_s^n|ds\nonumber\\
&+&C \int_t^T|Y_{s}^n-Y_s||Y_{s^+}^n-Y_s^n||\int_s^{s^+}Z_v^ndW_v||\dot{B}_s^n|ds\nonumber\\
&:=&B_{31,1}+B_{31,2}+B_{31,3}.
\end{eqnarray}
Furthermore,
\begin{eqnarray}\label{2.62}
B_{31,1}&\leq &C \int_t^T|\int_s^{s^+}f(Y^n_v,Z_v^n)dv||\int_s^{s^+}Z_v^ndW_v||\dot{B}_s^n|ds\nonumber\\
&+&C \int_t^T|\int_s^{s^+}g(Y^n_v)dB_v^n||\int_s^{s^+}Z_v^ndW_v||\dot{B}_s^n|ds\nonumber\\
&+&C \int_t^T|\int_s^{s^+}Z_v^ndW_v|^2|\dot{B}_s^n|ds\nonumber\\
&:=&B_{31,11}+B_{31,12}+B_{31,13}.
\end{eqnarray}
Now, interchanging the order of integration, we have
\begin{eqnarray}\label{2.63}
E[B_{31,11}]&\leq &C \frac{1}{2^n}\int_t^T(E[|\int_s^{s^+}(Z_v^n)^2dv|^2])^{\frac{1}{2}} (E[|\dot{B}_s^n|^2])^{\frac{1}{2}}ds\nonumber\\
&\leq& C \frac{1}{2^n}.
\end{eqnarray}
In view of (\ref{2.47}), (\ref{2.48}), we have
\begin{eqnarray}\label{2.64}
E[B_{31,12}]&\leq &C \frac{1}{2^n} \int_t^TE[\sup_s|\dot{B}_s^n|^2 |\int_s^{s^+}Z_v^ndW_v|]ds\nonumber\\
&\leq &C \frac{1}{2^n} (2^n)^2\int_t^T(E[\sup_{|r-v|\leq \frac{1}{2^n}} |B_r-B_v|^4])^{\frac{1}{2}} (E[|\int_s^{s^+}Z_v^ndW_v|^2])^{\frac{1}{2}} ds\nonumber\\
&\leq &C \frac{1}{2^n} (2^n)^2(\frac{1}{2^n})^{1-\delta}(\int_t^T E[\int_s^{s^+}(Z_v^n)^2dv]ds)^{\frac{1}{2}}\nonumber\\
&\leq &C \frac{1}{2^n} (2^n)^2(\frac{1}{2^n})^{1-\delta}(\frac{1}{2^n})^{\frac{1}{2}}(\int_t^T E[(Z_v^n)^2]dv)^{\frac{1}{2}}\nonumber\\
&\leq& C(\frac{1}{2^n})^{\frac{1}{2}-\delta}.
\end{eqnarray}
Using the fact that $|\dot{B}_s^n|$ is ${\cal F}_s$-measurable we have
\begin{eqnarray}\label{2.65}
E[B_{31,13}]&\leq &CE[\int_t^T(\int_s^{s^+}Z_v^ndW_v)^2|\dot{B}_s^n|ds] \nonumber\\
&=&CE[\int_t^T(\int_s^{s^+}(|\dot{B}_s^n|)^{\frac{1}{2}}Z_v^ndW_v)^2ds] \nonumber\\
&=&CE[\int_t^T\int_s^{s^+}(|\dot{B}_s^n|)(Z_v^n)^2dv ds] \nonumber\\
&\leq &C (2^n)E[(\sup_{|r-v|\leq \frac{1}{2^n}} |B_r-B_v|) \int_t^T\int_s^{s^+}((Z_v^n)^2dv ds]\nonumber\\
&\leq &C (2^n)(E[(\sup_{|r-v|\leq \frac{1}{2^n}} |B_r-B_v|^2])^{\frac{1}{2}} (E[(\int_t^T\int_s^{s^+}(Z_v^n)^2dv ds)^2])^{\frac{1}{2}} \nonumber\\
&\leq& C(\frac{1}{2^n})^{\frac{1}{2}-\delta},
\end{eqnarray}
where (\ref{2.47}),(\ref{2.48}) again were used. (\ref{2.63})- (\ref{2.65}) implies that
\begin{equation}\label{2.65-1}
  E[B_{31,1}]\leq C(\frac{1}{2^n})^{\frac{1}{2}-\delta}.
  \end{equation}
  For the term $B_{31,2}$ we have
\begin{eqnarray}\label{2.66}
B_{31,2}&\leq &C \int_t^T|\int_s^{s^+}f(Y_v,Z_v)dv||\int_s^{s^+}Z_v^ndW_v||\dot{B}_s^n|ds\nonumber\\
&+&C \int_t^T|\int_s^{s^+}g(Y_v)dB_v||\int_s^{s^+}Z_v^ndW_v||\dot{B}_s^n|ds\nonumber\\
&+&C \int_t^T|\int_s^{s^+}Z_vdW_v||\int_s^{s^+}Z_v^ndW_v||\dot{B}_s^n|ds\nonumber\\
&:=&B_{31,21}+B_{31,22}+B_{31,23}.
\end{eqnarray}
By a similar argument as for (\ref{2.63}), we have
\begin{equation}\label{2.67}
E[B_{31,21}]\leq C(\frac{1}{2^n}).
\end{equation}
As for $B_{31,22}$, we have
\begin{eqnarray}\label{2.68}
E[B_{31,22}]&\leq &C \int_t^T(E[|\int_s^{s^+}Z_u^n(|\dot{B}_s^n|)^{\frac{1}{2}}dW_u|^2])^{\frac{1}{2}}
(E[|\int_s^{s^+}g(Y_u)dB_u|^2])^{\frac{1}{2}}ds\nonumber\\
&\leq &C(\frac{1}{2^n})^{\frac{1}{2}} \int_t^T(E[\int_s^{s^+}(Z_u^n)^2|\dot{B}_s^n|du])^{\frac{1}{2}}ds\nonumber\\
\end{eqnarray}
From here following the same arguments as for (\ref{2.65}) we obtain
\begin{equation}\label{2.69}
E[B_{31,22}]\leq C(\frac{1}{2^n})^{\frac{1}{2}-delta}.
\end{equation}
The term $B_{31,23}$ is bounded as
\begin{eqnarray}
B_{31,23}&\leq &C \int_t^T(\int_s^{s^+}Z_v^ndW_v)^2|\dot{B}_s^n|ds\nonumber\\
&+ &C \int_t^T(\int_s^{s^+}Z_vdW_v)^2|\dot{B}_s^n|ds,
\end{eqnarray}
which, together with the same arguments as for (\ref{2.84}), yields
\begin{equation}\label{2.70}
E[B_{31,23}]\leq C(\frac{1}{2^n})^{\frac{1}{2}-\delta}.
\end{equation}
(\ref{2.67})- (\ref{2.70}) gives that
\begin{equation}\label{2.70-1}
  E[B_{31,2}]\leq C(\frac{1}{2^n})^{\frac{1}{2}-\delta}.
  \end{equation}
Finally we need to find an upper bound for $E[B_{31,3}]$. Notice that
\begin{eqnarray}\label{2.71}
B_{31,3}&\leq &C \int_t^T|Y_s^n-Y_s||\int_s^{s^+}f(Y^n_v,Z_v^n)dv||\int_s^{s^+}Z_v^ndW_v||\dot{B}_s^n|ds\nonumber\\
&+&C \int_t^T|Y_s^n-Y_s||\int_s^{s^+}g(Y^n_v)dB_v^n||\int_s^{s^+}Z_v^ndW_v||\dot{B}_s^n|ds\nonumber\\
&+&C \int_t^T|Y_s^n-Y_s||\int_s^{s^+}Z_v^ndW_v|^2|\dot{B}_s^n|ds\nonumber\\
&:=&B_{31,31}+B_{31,32}+B_{31,33}.
\end{eqnarray}
We have
\begin{eqnarray}\label{2.72}
E[B_{31,31}]&\leq &C \frac{1}{2^n} \int_t^T(E[|Y_s^n-Y_s|^2])^{\frac{1}{2}} (E[|\int_s^{s^+}(|\dot{B}_s^n|)^{\frac{1}{2}}Z_v^ndW_v|^2])^{\frac{1}{2}} ds\nonumber\\
&\leq &C \frac{1}{2^n} (\int_t^TE[|Y_s^n-Y_s|^2]ds)^{\frac{1}{2}} (\int_t^TE[|\int_s^{s^+}|\dot{B}_s^n|(Z_v^n)^2dv]ds)^{\frac{1}{2}}\nonumber\\
&\leq &C \frac{1}{2^n}(2^n)^{\frac{1}{2}}(E[(\sup_{|r-v|\leq \frac{1}{2^n}} |B_r-B_v|) \int_t^T\int_s^{s^+}(Z_v^n)^2dv]ds)^{\frac{1}{2}}\nonumber\\
&\leq &C (\frac{1}{2^n})^{\frac{1}{2}}(E[\sup_{|r-v|\leq \frac{1}{2^n}} |B_r-B_v|^2])^{\frac{1}{4}}(E[( \int_t^T\int_s^{s^+}(Z_v^n)^2dv ds)^2])^{\frac{1}{4}}\nonumber\\
&\leq &C (\frac{1}{2^n})^{\frac{1}{2}-\delta}(E[( \int_t^T(Z_v^n)^2dv)^2])^{\frac{1}{4}}\nonumber\\
&\leq &C (\frac{1}{2^n})^{\frac{1}{2}-\delta},
\end{eqnarray}
where the a priori bounds (\ref{2.1}), (\ref{2.1-1}) have been used.
Noticing that $|Y_s^n-Y_s|,|\dot{B}_s^n|ds$ are ${\cal F}_s$ measurable we have
\begin{eqnarray}\label{2.73}
E[B_{31,33}]&= &C \int_t^TE[|Y_s^n-Y_s||\int_s^{s^+}Z_v^ndW_v|^2|\dot{B}_s^n|]ds\nonumber\\
&=&C \int_t^TE[|Y_s^n-Y_s|(\int_s^{s^+}(Z_v^n)^2dv)|\dot{B}_s^n|]ds \nonumber\\
&\leq &C2^n E[ (\sup_{0\leq s\leq T}|Y_s^n-Y_s|)(\sup_{|r-v|\leq \frac{1}{2^n}} |B_r-B_v|^2])\int_t^T(\int_s^{s^+}(Z_v^n)^2dv)ds]\nonumber\\
&\leq&C2^n (E[\sup_{0\leq s\leq T}|Y_s^n-Y_s|^4])^{\frac{1}{4}}(E[\sup_{|r-v|\leq \frac{1}{2^n}} |B_r-B_v|^4])^{\frac{1}{4}}\nonumber\\
&&\quad\quad\times (E[(\int_t^T(\int_s^{s^+}(Z_v^n)^2dv)ds)^2])^{\frac{1}{2}}\nonumber\\
&\leq&C2^n (\frac{1}{2^n})^{\frac{1}{2}-\delta}(E[(\int_t^T(Z_v^n)^2
(\int_{v^-}^{v}ds)dv)^2])^{\frac{1}{2}}\nonumber\\
&\leq&C2^n (\frac{1}{2^n})^{\frac{1}{2}-\delta}\frac{1}{2^n}(E[(\int_t^T(Z_v^n)^2dv)^2])^{\frac{1}{2}}
\nonumber\\
&\leq&C(\frac{1}{2^n})^{\frac{1}{2}-\delta}.
\end{eqnarray}
As for the term $B_{31.32}$ we have
\begin{eqnarray}\label{2.74}
E[B_{31,32}]&\leq &C \int_t^TE[|Y_s^n-Y_s||\int_s^{s^+}Z_v^ndW_v|^2|\dot{B}_s^n|]ds\nonumber\\
&+ &C \int_t^TE[|Y_s^n-Y_s||\int_s^{s^+}g(Y_v^n)dB_v^n|^2|\dot{B}_s^n|]ds.
\end{eqnarray}
Furthermore, we have
\begin{eqnarray}\label{2.75}
&& \int_t^TE[|Y_s^n-Y_s||\int_s^{s^+}g(Y_v^n)dB_v^n|^2|\dot{B}_s^n|]ds\nonumber\\
&\leq& C(2^n)^3E[(\sup_{|r-v|\leq \frac{1}{2^n}} |B_r-B_v|^3)\int_t^T|Y_s^n-Y_s|(\int_s^{s^+}|g(Y_v^n)|dv)^2ds]\nonumber\\
&\leq& C(2^n)(E[(\sup_{|r-v|\leq \frac{1}{2^n}} |B_r-B_v|^6])^{\frac{1}{2}}(E[(\int_t^T|Y_s^n-Y_s|ds)^2])^{\frac{1}{2}}\nonumber\\
&\leq&C(\frac{1}{2^n})^{\frac{1}{2}-\delta}.
\end{eqnarray}
It follows now from (\ref{2.73}), (\ref{2.74}) and (\ref{2.75}) that
\begin{equation}\label{2.76}
E[B_{31,32}]\leq C(\frac{1}{2^n})^{\frac{1}{2}-\delta}.
\end{equation}
(\ref{2.72})- (\ref{2.76}) yields that
\begin{equation}\label{2.76-1}
  E[B_{31,3}]\leq C(\frac{1}{2^n})^{\frac{1}{2}-\delta}.
\end{equation}
It follows from (\ref{2.65-1}), (\ref{2.70-1}) and (\ref{2.76-1}) that
\begin{equation}\label{2.77}
  E[B_{31}]\leq C(\frac{1}{2^n})^{\frac{1}{2}-\delta}.
\end{equation}
This together with (\ref{2.60}) yields
\begin{equation}\label{2.78}
  E[B_{3}]\leq C(\frac{1}{2^n})^{\frac{1}{2}-\delta}.
\end{equation}
The lemma now follows from (\ref{2.44}), (\ref{2.58-1}) and (\ref{2.78}).
\vskip 0.4cm
\noindent {\bf Proof of Theorem 2.3}(continued).

  We are ready to complete the proof of Theorem 2.3. Taking expectation in (\ref{2.17})
  we obtain
\begin{eqnarray}\label{2.79}
&&E[(Y_t^n-Y_t)^2]+E\left [\int_t^T(Z_s^n-Z_s)^2ds\right ]\nonumber\\
&=& 2E\left[\int_t^T(Y_s^n-Y_s)(f(Y_s^n,Z_s^n)-f(Y_s,Z_s)) ds\right]\nonumber\\
&+&2E\left[\int_t^T(Y_s^n-Y_s)g(Y_s^n)dB_s^n\right]+E\left[\int_t^Tg^2(Y_s)ds\right]\nonumber\\
&-&E\left [\int_t^T(Y_s^n-Y_s)gg^{\prime}(Y_s)ds\right ]\nonumber\\
\end{eqnarray}
Taking into account of the  estimates in Lemma 2.4 and  Lemma 2.5
  we deduce from (\ref{2.79}) that
  \begin{eqnarray}\label{2.80}
&&E[(Y_t^n-Y_t)^2]+E\left [\int_t^T(Z_s^n-Z_s)^2ds\right ]\nonumber\\
&\leq &C(\frac{1}{2^n})^{\frac{1}{2}-\delta}+ CE\left[\int_t^T|Y_s^n-Y_s|^2 ds\right]\nonumber\\
&+& \frac{1}{2}E\left[\int_t^T|Z_s^n-Z_s|^2 ds\right]+E\left[\int_t^Tg^2(Y_s)ds\right]\nonumber\\
&-&E\left [\int_t^T(Y_s^n-Y_s)gg^{\prime}(Y_s)ds\right ]\nonumber\\
&+&E[\int_t^T(Y^n_{s^+}-Y_{s^+})gg^{\prime}(Y^n_{s^+})ds]\nonumber\\
&&-2E[\int_t^Tg(Y_{s^+})g(Y^n_{s^+})ds]+E[\int_t^Tg^2(Y^n_{s^+})ds].
\end{eqnarray}
To proceed with the proof, we claim that there is a constant $C$ such that
\begin{equation}\label{2.81}
E[\int_0^T|Y^n_{s^+}-Y_{s}^n|^2ds]\leq C(\frac{1}{2^n})^{1-\delta}.
\end{equation}
\begin{equation}\label{2.82}
E[\int_0^T|Y_{s^+}-Y_{s}|^2ds]\leq C(\frac{1}{2^n})^{1-\delta}.
\end{equation}
Let us prove (\ref{2.81}). The proof of (\ref{2.82}) is similar. Indeed, we have
\begin{eqnarray}\label{2.83}
&&E[\int_0^T|Y^n_{s^+}-Y_{s}^n|^2ds]\nonumber\\
&\leq &CE[\int_0^T|\int_{s}^{s^+}f(Y_u^n,Z_u^n)du|^2ds]\nonumber\\
&+&CE[\int_0^T|\int_{s}^{s^+}g(Y_u^n)dB_u^n|^2ds]\nonumber\\
&+&CE[\int_0^T|\int_{s}^{s^+}Z_u^ndW_u|^2ds]\nonumber\\
&\leq &C\frac{1}{2^n}+CE[\sup_{|r-v|\leq \frac{1}{2^n}}|B_r-B_v|^2]\nonumber\\
&+&CE[\int_0^T\int_{s}^{s^+}(Z_u^n)^2duds]\nonumber\\
&\leq &C(\frac{1}{2^n})^{1-\delta}+\sup_n (E[\int_0^T(Z_u^n)^2du])\frac{1}{2^n}\nonumber\\
&\leq& C(\frac{1}{2^n})^{1-\delta}.
\end{eqnarray}
By H\"o{}lder inequality it follows immediately from (\ref{2.81}) and (\ref{2.82}) that
\begin{equation}\label{2.84}
E[\int_0^T|Y^n_{s^+}-Y_{s}^n|ds]\leq C(\frac{1}{2^n})^{\frac{1}{2}-\delta}.
\end{equation}
\begin{equation}\label{2.85}
E[\int_0^T|Y_{s^+}-Y_{s}|ds]\leq C(\frac{1}{2^n})^{\frac{1}{2}-\delta}.
\end{equation}
Because of (\ref{2.81}) and (\ref{2.82}), we now replace $s^+$ by $s$ on the right side of
(\ref{2.80}) to obtain
\begin{eqnarray}\label{2.86}
&&E[(Y_t^n-Y_t)^2]+E\left [\int_t^T(Z_s^n-Z_s)^2ds\right ]\nonumber\\
&\leq &C(\frac{1}{2^n})^{\frac{1}{2}-\delta}+ CE\left[\int_t^T|Y_s^n-Y_s|^2 ds\right]\nonumber\\
&+& \frac{1}{2}E\left[\int_t^T|Z_s^n-Z_s|^2 ds\right]+E\left[\int_t^Tg^2(Y_s)ds\right]\nonumber\\
&-&E\left [\int_t^T(Y_s^n-Y_s)gg^{\prime}(Y_s)ds\right ]\nonumber\\
&+&E[\int_t^T(Y^n_{s}-Y_{s})gg^{\prime}(Y^n_{s})ds]\nonumber\\
&&-2E[\int_t^Tg(Y_{s})g(Y^n_{s})ds]+E[\int_t^Tg^2(Y^n_{s})ds].
\end{eqnarray}
Remark that the constant $C$ in front of $(\frac{1}{2^n})^{\frac{1}{2}-\delta}$ is different from that in (\ref{2.80}). Completing the square in (\ref{2.86}) we get that
\begin{eqnarray}\label{2.87}
&&E[(Y_t^n-Y_t)^2]+E\left [\int_t^T(Z_s^n-Z_s)^2ds\right ]\nonumber\\
&\leq &C(\frac{1}{2^n})^{\frac{1}{2}-\delta}+ CE\left[\int_t^T|Y_s^n-Y_s|^2 ds\right]\nonumber\\
&+& \frac{1}{2}E\left[\int_t^T|Z_s^n-Z_s|^2 ds\right]+E\left[\int_t^T(g(Y_s)-g(Y_s^n))^2ds\right]\nonumber\\
&+&E[\int_t^T(Y^n_{s}-Y_{s})(gg^{\prime}(Y^n_{s})-gg^{\prime}(Y_s))ds].
\end{eqnarray}
By the Lipschitz continuity of $gg^{\prime}$ and $g$, it follows from (\ref{2.87}) that
\begin{eqnarray}\label{2.88}
&&E[(Y_t^n-Y_t)^2]+\frac{1}{2} E\left [\int_t^T(Z_s^n-Z_s)^2ds\right ]\nonumber\\
&\leq &C(\frac{1}{2^n})^{\frac{1}{2}-\delta}+ CE\left[\int_t^T|Y_s^n-Y_s|^2 ds\right].
\end{eqnarray}
Application of the Gronwall's inequality completes the proof of Theorem 2.3.

\end{document}